\documentclass[smallextended]{svjour3}      

\smartqed  

\usepackage{graphicx}
\usepackage{amssymb}
\usepackage{amsmath}
\usepackage{color}
\usepackage{epstopdf}
\usepackage{enumitem}

\begin{document}

\title{An estimation of distribution algorithm for the computation of innovation estimators of dif\mbox{}fusion processes}

\author{Zochil Gonz\'alez Arenas         \and
        Juan Carlos Jim\'enez            \and 
        Li-Vang Lozada-Chang             \and
        Roberto Santana
}

\institute{Z.G. Arenas \at
              Institute of Mathematics and Statistics, University of the State of Rio de Janeiro, Brazil\\
              Tel.: +55-21-993063131\\
              \email{zochilita@gmail.com; zochil@uerj.br}          
           \and
           J.C. Jimenez \at
              Institute for Cybernetics, Mathematics and Physics, Havana, Cuba
           \and
           L. Lozada-Chang \at
              Faculty of Mathematics and Computation, University of Havana, Cuba
           \and
           R. Santana \at
             Dept. of Computer Science and  Artificial Intelligence, University of the  Basque Country, Spain
}

\date{Received: date / Accepted: date}

\maketitle

\begin{abstract}
Estimation of Distribution Algorithms (EDAs) and Innovation Method
are recognized methods for solving global optimization problems and
for the estimation of parameters in diffusion processes,
respectively. Well known is also that the quality of the Innovation
Estimator strongly depends on an adequate selection of the initial
value for the parameters when a local optimization algorithm is used
in its computation. Alternatively, in this paper, we study the
feasibility of a specific EDA - a continuous version of the
Univariate Marginal Distribution Algorithm (UMDAc) - for the
computation of the Innovation Estimators. Numerical experiments are performed
for two different models with a high level of complexity.
The numerical simulations
show that the considered global optimization algorithms
substantially improves the effectiveness of the Innovation
Estimators for different types of diffusion processes with complex
nonlinear and stochastic dynamics.
\keywords{Local Linear Approximation method \and Estimation of Distribution Algorithms \and Parameter estimation \and Numerical simulations}

\end{abstract}

\section{Introduction}
\label{intro}

Dif\mbox{}fusion processes defined through Stochastic Differential Equations (SDEs) have became an important mathematical tool for
describing the dynamics of several phenomena, e.g., the dynamics of assets prices in the market, the fire of neurons, etc. In many applications,
the statistical inference of dif\mbox{}fusion processes is of great importance for model building and model
selection. This inference problem consists in the estimation of the unknown parameters and unobserved components of the dif\mbox{}fusion
process given a set of discrete and noisy observations of some of its components. In this context, a variety of inference methods
(see, e.g., \cite{Jim06a} ) have been considered and, among them, the Innovation Estimators have been shown very useful in
applications (see, e.g., \cite{Vald99,Oza2000,Riera07}).

The computation of the Innovation Estimators involves the minimization of an objective or fitness function, which is a non-quadratic
function of the parameters in most of situations. This optimization process is usually carried out by means of local optimization
algorithms \cite{Vald99,Oza2000,Riera07}. In general, innovation estimators computed in that way strongly depend on the quality of
the parameter's initial value used by the local optimization algorithm. Therefore, great expertise of the users is needed to
reach satisfactory estimates.

To overcome such type of difficulty, instead of a local optimization algorithm, global optimization methods as the Estimation of
Distribution Algorithms (EDAs) \cite{Larranaga_et_al:2012,Larr2002,Lozano_et_al:2005,Muhl1996} are usually considered. In particular, 
EDAs comprise a group of stochastic optimization heuristics which base the search of
an optimal solution on a population of individuals. In the population, each one of the individuals represents a solution to the
considered optimization problem. Individuals are evaluated and a subset of them is selected according to the quality of their 
objective function values. EDAs use probabilistic models of the solutions to extract relevant information about the set of selected
solutions. The probabilistic model is used to sample new individuals. In this way, the algorithm evolves in successive
generations towards the more promising regions of the search space until a stopping criterion is satisfied. 
A pseudocode of a general EDA will be later presented.

EDAs can be classified according to the type of solution representation and to the way that
the learning of the probability model is accomplished. The representations are discrete or continuous. Moreover, regarding the way of
learning, EDAs can be classified in two classes. One class groups the algorithms that make a parametric learning of the probabilities,
and the other one comprises those algorithms where structural as well as parametric learning of the model is done.

In this paper, for the computation of the Innovation Estimators of the unknown parameters of dif\mbox{}fusion processes, we focus on
a class of EDAs with continuous representation and parametric learning. As probabilistic model, we use the univariate marginals
estimated from the selected population, which defines the so called Univariate Marginal Distribution Algorithm in continuous domain
(UMDAc) \cite{Muhl98}. This global optimization strategy is also combined with a local optimization algorithm and tested in numerical
simulations.

The paper is organized as follows. In Section~\ref{sec:Estimation}, the estimation
problem is clearly defined, and the essentials on the Innovation
Method and EDAs are briefly presented. Section~\ref{sec:UMDAcInnovEst} focused on the
application of the UMDAc  to the parameter estimation of SDEs and
the resulting algorithms are summarized. The performance of the
proposed algorithms is illustrated in Section~\ref{sec:EXPER} in the estimation of
two types of diffusion processes with complex nonlinear and
stochastic dynamics. Finally, in Section~\ref{sec:CONCLU}, we present the conclusions of the paper
and discuss some possible lines of future work.


\section{Notations and Preliminaries}  
\label{sec:Estimation}

\bigskip Let the continuous-discrete state space model be defined by the continuous state equation
\begin{equation}
dx(t) = f(t,x(t);\alpha )dt+\sum\limits_{i=1}^{m} g_i(t,x(t);\alpha
)dw^i(t),\text{ for }t\geq t_0\in \Re.\label{SS1}
\end{equation}
and the discrete observation equation
\begin{equation}
z_{t_k}=h_0(t_k,x(t_k))+e_{t_k},\ \text{ for } k=0,1,..,N\in
\mathbb{N},  \label{SS3}
\end{equation}
where $x(t)\in \Re ^d$ is the state vector at the instant of time $t$, $z_{t_k}\in \Re ^r$ is the observation vector at the instant of
time $t_{k}$, $\alpha $ is a set of $p$ parameters, $w$ is an \textit{m}-dimensional standard Wiener process,
$\{e_{t_k}:e_{t_k}\thicksim \mathcal{N}(0,\Sigma)\}$, $k=\{0,..,N\}$, is a sequence of i.i.d. random vectors independent
of $w$, $f,\ g_i:\ \Re \times \Re ^d\rightarrow \Re ^d$ are vector functions and $h_0:\ \Re \times \Re ^d\rightarrow \Re ^r$ is also a
vector function. Here, the time discretization $\left( t\right)_{N}=\left\{ t_{k}:k=0,1,\cdots ,N\right\} $ is assumed to be
increasing, i.e., $t_{k-1}<t_{k}$ for all $k=1,..,N$.

Let $x_{t/t_{k}}=\mathbf{E}(x(t)/Z_{t_{k}};\alpha)$ and 
$P_{t/t_k}=\mathbf{E}((x(t)-x_{t/t_k}) (x(t)-x_{t/t_{k}})^{\intercal}/Z_{t_{k}};\alpha )$ for all $t\geq t_{k}$, where $\mathbf{E}(./.)$
denotes conditional expectation and $Z_{t_{k}}=\{z_{t_{j}}:$
$t_{j}\leq t_{k}$ and $t_{j}\in \left(t\right) _{N}\}$ is a sequence of observations. In the case that $t>t_{k}$, $x_{t/t_{k}}$ and
$P_{t/t_{k}}$ are called prediction and prediction variance,
respectively. If $t=t_{k}$, $x_{t_{k}/t_{k}}$ and $P_{t_{k}/t_{k}}$
are called filter and filter variance. Let $\vartheta
=\{x_{t_{0}/t_{0}},P_{t_{0}/t_{0}},\alpha \}$.

The inference problem here consists in the estimation of
$\vartheta$, $x_{t_{k}/t_{k}}$ and $P_{t_{k}/t_{k}}$ given the time
series $Z_{t_{N}}$ with $N$ observations $z_{t_{k}}$. Usually, the
inference problem for differential equations is carried out in two
steps: 1) the estimation of the unknown parameters $\alpha$, and 2)
the estimation of the unobserved component $x$ at each $t_k$, i.e.,
all the $x_{t_{k}/t_{k}}$ and $P_{t_{k}/t_{k}}$. In a model like
(\ref{SS1})-(\ref{SS3}) with $\alpha$ given, the problem of
estimating the state $x(t_k)$ from the observations $z_{t_k}$ is
known as the nonlinear continuous-discrete f\mbox{}iltering problem.
In practical applications, $\alpha$ is commonly unknown. For this
reason, the estimation of these parameters plays a central role.

\section[Innovation Estimators]{Innovation Estimators}  \label{sec:InnovEstimators}

The Innovation Estimator $\widehat{\alpha}$ of the parameter
$\alpha$ in the model (\ref{SS1})-(\ref{SS3}) is defined as
\cite{Jim06b}:
\begin{equation}
\widehat{\alpha}=\arg {\underset{\alpha}{\min} \left\{ q(\alpha)
\right\}}, \label{eq:InnovEstim}
\end{equation}
where
\[q(\alpha) = \left\{ N\ln (2\pi )\mathbf{+}\sum\limits_{k=1}^{N}\ln (\det (\Sigma_{t_k/t_{k-1}}^{\mathbf{\nu }}))+\nu_{t_k}^{\top }(\Sigma_{t_k/t_{k-1}}^{\mathbf{\nu }})^{-1}\nu_{t_k}\right\},\]
$\nu_{t_k}=z_{t_k}-\mathbf{E}(h(t_{k},x(t_k))/Z_{t_{k}};\alpha)$ and $\Sigma_{t_k/t_{k-1}}^{\mathbf{\nu }}$ are, 
respectively, the discrete
time innovation process and its variance.

Since the exact computation of $\nu_{t_k}$ and
$\Sigma_{t_k/t_{k-1}}$ is only possible for a few particular models,
in general, it is necessary to use approximate formulas. The
approximations $\widetilde{\nu}_{t_{k}}$ and
$\widetilde{\Sigma}_{t_k/t_{k-1}}^{\mathbf{\nu }}$ are recursively
computed with the following Local Linearization filtering
algorithm~\cite{Jim03}, for $k=0,..,N-1$:

\begin{enumerate}
\item Prediction
\begin{equation*}
y_{t_{k+1}/t_k}=y_{t_k/t_k}+\int\limits_0^{t_{k+1}-t_k}(A_k y_{(t_k+t)/t_k}+a_k(t))dt
\label{LLF1}
\end{equation*}
\begin{eqnarray*}
Q_{t_{k+1}/t_k} & = & Q_{t_k/t_k}+\int\limits_0^{t_{k+1}-t_k}\{A_kQ_{(t_k+t)/t_k}+
Q_{(t_k+t)/t_k}A_k^{\intercal}  \\
&&+\ \sum\limits_{i=1}^{m}B_{i,k}(Q_{(t_k+t)/t_k}+y_{(t_k+t)/t_k}y_{(t_k+t)/t_k}^{\intercal})
B_{i,k}^{\intercal } \\
&&+\ \sum\limits_{i=1}^{m}B_{i,k}y_{(t_{k}+t)/t_{k}}b_{i,k}^{\intercal}(t)+ b_{i,k}(t)y_{(t_{k}+t)/t_{k}}^{\intercal }B_{i,k}^{\intercal } \\ 
&&+\ \sum\limits_{i=1}^{m}b_{i,k}(t)b_{i,k}^{\intercal}(t)\}dt,
\label{LLF2}
\end{eqnarray*}

\item Innovation
\begin{equation*}
\widetilde{\nu}_{t_{k+1}}=z_{t_{k+1}}-h_0(y_{t_{k+1}/t_k}), 
\end{equation*}
\begin{equation*}
\widetilde{\Sigma}_{t_{k+1}/t_k}^{\nu}=C_kQ_{t_{k+1}/t_k}C^{\intercal}_k+\Sigma,
\end{equation*}

\item F\mbox{}ilter
\begin{equation*}
y_{t_{k+1}/t_{k+1}}=y_{t_{k+1}/t_k}+K_{t_{k+1}}\widetilde{\nu}_{t_{k+1}},
\label{LLF3}
\end{equation*}
\begin{equation*}
Q_{t_{k+1}/t_{k+1}}=Q_{t_{k+1}/t_k}-K_{t_{k+1}}C_kQ_{t_{k+1}/t_k},
\label{LLF4}
\end{equation*}
where
$K_{t_{k+1}}=Q_{t_{k+1}/t_k}C_k^{\intercal} \left(\widetilde{\Sigma}_{t_{k+1}}^{\nu}\right)^{-1}$ is the f\mbox{}ilter gain.
\label{problem}
\end{enumerate}

\vspace{.4cm}
In the above algorithm the remaining notations are:
\[ A_k=J_f(t_k,y_{t_k/t_k}), \qquad B_{i,k}=J_{g_i}(t_k,y_{t_k/t_k}),\quad C_k=J_{h_0}(t_k,y_{t_k/t_k}) \]

\begin{equation*}
a_k(t)=f(t_k,y_{t_k/t_k})-J_f(t_k,y_{t_k/t_k})y_{t_k/t_k}+
\frac{\partial f(t_k,y_{t_k/t_k})}{\partial s}(t-t_k),
\end{equation*}
and
\begin{equation*}
b_{i,k}(t)=g_i(t_k,y_{t_k/t_k})-J_{g_i}(t_k,y_{t_k/t_k})y_{t_k/t_k}+
\frac{\partial g_i(t_k,y_{t_k/t_k})}{\partial s}(t-t_k),
\end{equation*}
where $J_v$ indicates the Jacobian matrix of the vector function
$v$. The algorithm starts with
\[y_{t_0/t_0}=x_{t_0/t_0},\quad Q_{t_0/t_0}=P_{t_0/t_0}\]
for each given value of $\alpha$.

Finally, the estimates $y_{t_k/t_k}$ and $Q_{t_k/t_k}$ obtained from
the above filtering algorithm with $\alpha =\widehat{\alpha}$ give
an approximation of the mean $x_{t_{k}/t_{k}}$ and variance
$P_{t_{k}/t_{k}}$ of the states $x(t_k)$ in the
model~(\ref{SS1})-(\ref{SS3}).

Explicit formulas for the solution of the differential equations
involved in the prediction step can be found in~\cite{Jim03} or in~\cite{Jim15}.

The minimization of the function $q(\alpha)$ with respect to
$\alpha$ is a major difficulty in the computation of the Innovation
Estimators, as it can be observed from the
expression~(\ref{eq:InnovEstim}). In this situation, the
optimization process acquires great importance for estimating the
parameters of the model.

Among the difficulties of this optimization problem are the
non-quadratic dependence of the fitness or objective function
$q(\alpha)$ with respect to the parameters $\alpha$ due to the
nonlinearity of the model (\ref{SS1})-(\ref{SS3}) and the presence
of parameters $\alpha$ in the highly nonlinear term corresponding to
the innovation variance $\Sigma^\nu$. An additional difficulty is
the impossibility of using local optimization methods of high
convergence order since calculating the gradient or the Hessian of
the objective function with regard to the parameters is not
possible.

\subsection[UMDAc]{Univariate Marginal Distribution Algorithm} \label{sec:UMDAc}

A pseudocode of a general Estimation of Distribution Algorithm for solving optimization problems is
shown in Algorithm~\ref{alg:EDA}.

\begin{table}[ht]
  \caption{ Pseudocode of a general Estimation of Distribution Algorithm (EDA)}
  \label{alg:EDA}       
  \normalsize
  \begin{tabular}{p{14cm}}
   \hline\noalign{\smallskip}
  \hspace{5cm} Algorithm \ref{alg:EDA}: {\bf EDA}
  \smallskip \\
  \hline\hline
     \begin{enumerate}[leftmargin=0em, leftmargin=*, label=\arabic*, itemsep=1pt, parsep=1pt, font=\it]
       \item \hspace{.2cm} Set $t\Leftarrow 1$. Generate $M$ points randomly from the search space (Initial population  $D_0$).
       \item \hspace{.3cm}{\bf do \{ }
       \item   \hspace{.8cm} {Evaluate the fitness function at the $M$ points.}
       \item   \hspace{.8cm} {Select a set of points $D^{Se}_{t-1}$ according to a selection method.}
       \item   \hspace{.8cm} {Compute a probabilistic model from $D^{Se}_{t-1}$.}
       \item   \hspace{.8cm} {Sample $M$ new points according to the probability distribution previously learnt.}
       \item   \hspace{.8cm} {$t \Leftarrow t+1$}
       \item  \hspace{.7cm} {\bf \} until } Termination criteria are met 
      \end{enumerate}\\
\hline
\end{tabular}
\end{table}

In particular, the Univariate Marginal Distribution Algorithm (UMDA) \cite{Muhl98} follows Algorithm~\ref{alg:EDA}, but uses
univariate marginals calculated from the selected population as the probabilistic model.  Theoretical studies, as well as a variety of
applications of UMDA have been reported in \cite{Liv11,Hash2011,Muhl1996,San07}. On the other hand,
applications of EDAs to problems with continuous representation comprise the use of UMDA based on Gaussian models
\cite{Beng2002,Larr1999,Seb98}, Gaussian networks \cite{Beng2002,Lar2015,Larr1999}, Marginal Histogram in Continuous
Domains~\cite{Tsut2001}, mixtures of Gaussian distributions \cite{Bos2000a}, and Voronoi based EDAs \cite{Oka2004}. In its
general form, the UMDA for continuous domains (UMDAc)~\cite{Larr1999} is not restricted to the use of Gaussian
distributions, the density function which better fits the optimal
solutions is statistically determined for each generation. For more
details on EDAs for continuous domains \cite{Bosman_and_Grahl:2008,Bosman_and_Thierens:2006,Pelik99c} can
be consulted.

\section[UMDAc-based Innovation Estimators]{UMDAc-based Innovation Estimators for Dif\mbox{}fusions} 
\label{sec:UMDAcInnovEst}

Let $\Omega = \Re^p$ be the \emph{search space}, $X$ a random element of $\Omega$, and $X_i$ its $i$-th component.
$F(X_i;\theta_i)$ and $F(X;\theta)$ will denote, respectively, the density function of $X_i$ and the joint density function of $X$
depending of the sets of parameters $\theta_i$ and $\theta = \left\{\theta_1,..,\theta_p \right\}$.

In addition, $D_l$ and $D^{Se}_{l}$ will denote the population at the $l$-th generation and the selected population at the $l$-th
generation from which the joint probability distribution of $X$ is learnt. $F_{l}\left(X;\widehat\theta^l\right)$ denotes the joint
density function estimated in each generation $l$.

The pseudocode for learning the joint density function by the UMDAc at each generation $l$ is shown in Algorithm~\ref{alg:UMDAc}.

\begin{table}[ht]
  \caption{Calculation of the probabilistic model by the Univariate Marginal Distribution Algorithm for continuous domains}
  \label{alg:UMDAc}     
  \normalsize
  \begin{tabular}{p{15cm}}
   \hline\noalign{\smallskip}
  {\hspace{4cm} Algorithm \ref{alg:UMDAc}: {\bf  UMDAc: probabilistic model}}
  \smallskip\\
  \hline\hline
     \begin{enumerate}[leftmargin=0em, leftmargin=*, label=\arabic*, itemsep=0pt, parsep=0pt, font=\it]
       \item   \hspace{.1cm} {\bf for} {$i := 1$} {\bf to} {$p$} {\bf \{ } 
       \item   \hspace{.5cm} \emph{Select} via hypothesis test the density function
       $F_l\left(X_{i};\theta^{\ l}_{i}\right)$ which better fits the population $D^{Se}_{l-1}$ projected on $X_i$
       \item   \hspace{.5cm}  \emph{Obtain} the maximum likelihood estimators $\widehat\theta^{\ l}_{i}$ for $\theta^l_{i}$ 
       \item[]   \hspace{.3cm} \textbf{\} }
      \end{enumerate}
      The learnt joint density function is expressed as 
      $F_{l}\left(X;\theta^l\right) = \prod^{p}_{i=1}F_{l}\left(X_{i};\widehat\theta^{\ l}_{i}\right)$\\
  \hline
  \end{tabular}
\end{table}

We will use the UMDAc with Gaussian distributions. By considering independence among the components of $X$ we assume that the joint
density function $F(X)$ is normal $p$-dimensional and could be factorized as the product of the normal univariate marginal
densities of each component $X_i$. Univariate marginal densities for each $X_i$ will be calculated as
\begin{equation}
 F(X_i) = \frac{1}{\sqrt{2\pi} \sigma_i}e^{ -\frac{(X_i- \mu_i)^2}{2 \sigma_i^2}}.
\label{eq:GAUSS}
\end{equation}

In our optimization problem each component $X_i$ belongs to a bounded interval $[a_i,b_i]$. A first population $D_0$ of $M$ points
is generated according a uniform distribution. In every generation $l \ (l>1)$ with $M$ points $X^{j}$, for every component $X_i$ we
performed the estimation of the mean $\mu_i$ and the standard deviation $\sigma_i$ by their sampling estimates,
 \[\hat{\mu_{i}} = \frac{1}{M}\sum^{M}_{j=1}X^{j}_{i}, ~~~
 \text{    and     } ~~~ \hat{\sigma_{i}} = \sqrt{\frac{1}{M}\sum^{M}_{j=1}\left(X^{j}_{i}-\hat{\mu_{i}}\right)^{2}}. \]

By considering elitism, the $\varepsilon$ points $X$ with the best solutions (minimum values of $q(\alpha)$) of the previous generation are
kept in the new population. New $(M-\varepsilon)$ points $X$ are randomly generated from the normal distribution with the estimated
parameters but retaining the i-th component inside the interval $[a_i,b_i]$. The new $X$ are added into the original population
replacing the non-selected ones. The process is repeated until a stop condition is met. Different stopping condition can be considered,
as a fixed number of generations or non-improvement criterion after a given number of generations, for example.
We considered a fixed number of generations, justified by some preliminary experiments and cost-benefit analysis.

In addition, we consider the strategy of combining the global and local optimization techniques in which the estimated values of the
parameters obtained by EDA are used as initial values of the local technique. The application of EDAs together with local optimization
techniques has been reported to notably improve the quality of the solutions for problems from different
domains \cite{Muehlenbein_and_Mahnig:2002,Pelikan:2005,Zhang_et_al:2003}.

For comparison purposes, we used the MATLAB function \emph{fmincon} as a Local Optimization Algorithm (LOA). This function searches for
the minimum of a nonlinear multivariate function with constraints, and needs an initial estimate for this task.

In what follows we summarize with a pseudocode the three estimation algorithms considered in the paper for computing the innovation
estimator (\ref{eq:InnovEstim}) of the parameters $\alpha$ in the model (\ref{SS1})-(\ref{SS3}).

\begin{table}[ht]
  \caption{Algorithm for global optimization}
  \label{alg:Global}     
  \normalsize
  \begin{tabular}{p{15cm}}
   \hline\noalign{\smallskip}
  {\hspace{4cm} Algorithm \ref{alg:Global}: {\bf UMDAc - Innovation estimator}}
  \smallskip\\
  \hline\hline
     \begin{enumerate}[leftmargin=0em, leftmargin=*, label=\arabic*, itemsep=1pt, parsep=1pt, font=\it]
       \item \hspace{.2cm} Set $t\Leftarrow 1$. Generate $M$ points $X$ using a uniform distribution (Initial population $D_0$ for $\alpha$).
       \item \hspace{.3cm}{\bf do \{ }
       \item   \hspace{.8cm} {Evaluate the fitness function $q(\alpha)$ at the $M$ points $X$.}
       \item   \hspace{.8cm} {Select a set of points $D^{Se}_{t-1}$ according to a selection method.}
       \item   \hspace{.8cm} {Compute a probabilistic model from $D^{Se}_{t-1}$ applying Algorithm~\ref{alg:UMDAc} with $F(X)$ in Eq.~\eqref{eq:GAUSS}.}
       \item   \hspace{.8cm} {Consider elitism: keeping the $\varepsilon$ points $X$ with the best solutions in the new population.}
       \item   \hspace{.8cm} {Sample $M - \varepsilon$ new points $X$ according to the probability distribution previously learnt.}
       \item   \hspace{.8cm} {$t \Leftarrow t+1$}
       \item  \hspace{.7cm} {\bf \} until } Termination criteria are met
       \item \hspace{.2cm} $\hat{\alpha}$ is the $X$ with the minimum value of $q(\alpha)$ in the last generation.
      \end{enumerate}\\
\hline
\end{tabular}
\end{table}

\begin{table}[ht]
  \caption{Global optimization with Local optimization refinement}
  \label{alg:GlobalLocal}    
  \normalsize
  \begin{tabular}{p{15cm}}
   \hline\noalign{\smallskip}
  {\hspace{3cm} Algorithm \ref{alg:GlobalLocal}: {\bf Refined UMDAc - Innovation estimator}}
  \smallskip\\
  \hline\hline
     \begin{enumerate}[leftmargin=0em, leftmargin=*, label=\arabic*, itemsep=1pt, parsep=1pt, font=\it]
       \item  \hspace{.2cm} Use \emph{Algorithm}~\ref{alg:Global} to find an estimate $\tilde{\alpha}$ of $\alpha$.
       \item  \hspace{.2cm}  Set $\alpha_0 \Leftarrow \tilde{\alpha}$.
       \item  \hspace{.2cm} Compute a new estimation $\hat{\alpha}$ of $\alpha$ with the MATLAB function \emph{fmincon} starting at $\alpha_0$.
     \end{enumerate}\\
\hline
\end{tabular}
\end{table}

\begin{table}[ht]
\caption{Local optimization strategy (LOA)}
\label{alg:Local}     
  \normalsize
  \begin{tabular}{p{15cm}}
   \hline\noalign{\smallskip}
  {\hspace{3cm} Algorithm \ref{alg:Local}: {\bf LOA - Innovation estimator}}
  \smallskip\\
  \hline\hline
     \begin{enumerate}[leftmargin=0em, leftmargin=*, label=\arabic*, itemsep=1pt, parsep=1pt, font=\it]
       \item \hspace{.2cm} Define $\alpha_0$ (Initial value for $\alpha$).
       \item \hspace{.2cm} Compute an estimation $\hat{\alpha}$ of $\alpha$ with the MATLAB function \emph{fmincon} starting at $\alpha_0$.
     \end{enumerate}\\
\hline
\end{tabular}
\end{table}

\section[Numerical Experiments and Results]{Numerical Experiments and Results} \label{sec:EXPER}

The objective of our experiments is to determine whether the use of UMDAc can improve the computation of Innovation Estimators for
unknown parameters of discrete observed diffusion processes. To do so, we compare the Algorithms \ref{alg:Global}, \ref{alg:GlobalLocal} and \ref{alg:Local} described above in the
search of the solution to the optimization problem described by Eq.~\eqref{problem}.

With these three algorithms, we carried out $100$ estimations of the parameters $\alpha $ for a given time series $Z_{t_{N}}$ of $N$
observations $z_{t_{k}}$ of two different types of diffusion processes. In each estimation, the initial population for
 Algorithms \ref{alg:Global} and \ref{alg:GlobalLocal}, and the initial value for the Algorithm \ref{alg:Local}
were randomly selected from a predefined set of possible values for each parameter. The population size $M$ for the UMDAc was decided in
agreement with the approach suggested in~\cite{Muhl01}, i.e., $M$ was set equal to $20$ times the number of parameters $p$ to be
estimated. In order to have a quicker convergence to the optimum, truncation selection was used as the selection method of
choice~\cite{Muhl01}. An empirical rule previously used in Factorized Distribution Algorithm (FDA) \cite{Muhl99a} 
to specify the truncation threshold $\tau$ locates it between $0.125$ and $0.4$. We fixed $\tau = 0.3$ since this 
choice of $\tau$ has been already reported in previous EDA applications in continuous domains~\cite{Cho04}. Elitism was implemented with $5\%$
of the population. A fixed number of generations was the stopping condition used although we also kept in mind a superior bound for
the value of the objective function. All these mentioned values were selected after conducting a set of preliminary experiments.

For the following two examples of diffusion processes, the
realization $Z_{t_{N}}$ of the model was computed by using the known
order 1 strong Local Linearization scheme for SDEs (see, e.g.,
\cite{Jim12}) on a time discretization finer than the observation
times $t_{0}$,...,$t_{N}$ and then subsampling the approximate
solution of $x$ at the time instants $t_{0}$,...,$t_{N}$.

\subsection[Additive noise example]{Stif\mbox{}f state equation with additive noise} \label{sec:stiff}

The modeling of stochastic behavior of neuronal populations has been
of great interest for several years. The stochastic Fitzhugh-Nagumo
equation has widely been used in analytical and simulation studies
of neuronal models~\cite{Tuck98,Tuck2003}. The biggest complexity to
deal with this model is due to the stiffness of the nonlinear state
equations~(\ref{eq:Kon1})-(\ref{eq:Kon2}), the components of the
state variable change at very different rates. This behavior results
in great difficulty for its numerical solution.

Let us consider the stochastic Fitzhugh-Nagumo model defined by the
continuous nonlinear state equations
\begin{eqnarray}
dx_1 &=& 100(x_1-\frac{x_1^3}{3}-x_2)dt \label{eq:Kon1} \\
dx_2 &=& \alpha_1 +\alpha_2 x_{1}dt+\alpha_3 dw_{2} \label{eq:Kon2}
\end{eqnarray}
and the discrete observation equation
\begin{equation}
    z_{t_k}=x(t_k)+e_{t_k}
\label{eq:KonObs}
\end{equation}
where $x(t) \in \Re^2$, $z_{t_k}\in \Re^2$, $\alpha_1 = 1$,
$\alpha_2 = 1$, $\alpha_3 = 0.1$, $x(0) = (x_1(0),x_2(0)) = (-0.9323, -0.6732)$
and $\{e_{t_{k}}:e_{t_{k}}\thicksim \mathcal{N}(0,10^{-6})\}$.

\begin{figure}[htb]
\centering
\includegraphics[width = 12cm, height = 10cm]{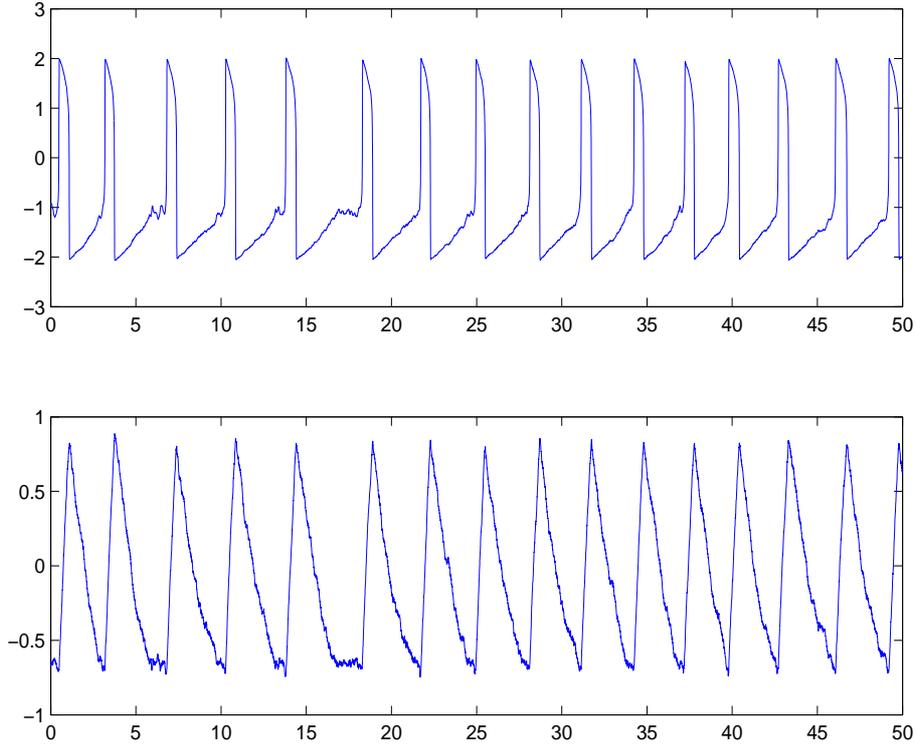}
\caption{ Realization of
the solution of the state space model (\ref{eq:Kon1})-(\ref{eq:Kon2}). Top:$\{x_1(t_j): j=1,\ldots,T\}.$
Bottom:$\{x_2(t_j): j=1,\ldots,T\}.$}
\label{fig:KonEq}
\end{figure}

A realization $\{x(t_j)=(x_1(t_j),x_2(t_j)): j=1,\ldots,T\}$ of the
solution of the state space model~(\ref{eq:Kon1})-(\ref{eq:Kon2}) is
shown in Figure~\ref{fig:KonEq} for instants of time $t_j=jh$,
with $h=0.0005$ and $T=50$. The variations of the state variable
which determines the stif\mbox{}f characteristics of the model are
observed.

Given the state space model~(\ref{eq:Kon1})-(\ref{eq:KonObs}) and
the single realization $\{z_{t_k}:k=0,\ldots ,N\}$ of the equation
of observations~(\ref{eq:KonObs}) with $t_k=k\Delta$, $\Delta =0.5$
and $N=500$, $100$ estimates of the parameter set $\alpha =
(\alpha_1,\alpha_2,\alpha_3)$ were carried out. Each initial
estimate of $\alpha$ was uniformly generated inside of their
definition intervals
\[\alpha_1 \in [0,5], \quad \alpha_2 \in [0,5], \quad
\alpha_3 \in [0,1].\]

Due to the complexity of this type of dif\mbox{}fusion process, it was not possible to carry out the estimation of the set of parameters
$\alpha$ using the Algorithm \ref{alg:Local}. To obtain a solution with the Matlab function \emph{fmincon} of Algorithm \ref{alg:Local} is essential to have a very
good initial estimation for the parameters. Otherwise, the algorithm does not converge. For this reason, the results of the estimation
with Algorithm \ref{alg:Local} is not reported in this example.

\begin{table}[ht]
    \centering
     \caption{\small Parameter estimation achieved by UMDAc for the state space model (\ref{eq:Kon1})-(\ref{eq:Kon2}).
      The interval is defined by the minimum and maximum values of the estimated parameters produced in $100$ executions of
      Algorithm~\ref{alg:Global}.}
      \normalsize
    \begin{tabular}{p{2cm} p{3cm} p{3cm}}
	\hline\hline\noalign{\smallskip}
        \centering Parameter & \centering{ Real value} & \hspace{.2cm} Estimated  \\ [.5ex]
        \hline\noalign{\smallskip}
        \centering $\alpha_1$ & \centering 1 & [0.9181,1.4887] \\ [.5ex]
        \centering $\alpha_2$ & \centering 1 & [0.9699,1.3701] \\ [.5ex]
        \centering $\alpha_3$ & \centering 0.1 & [0.1018,0.1458] \\ [1ex]
        \hline
    \end{tabular}
  \label{table:UMDAc}
\end{table}

\begin{figure}[ht]
\centering
\includegraphics[width = 12cm, height = 10cm]{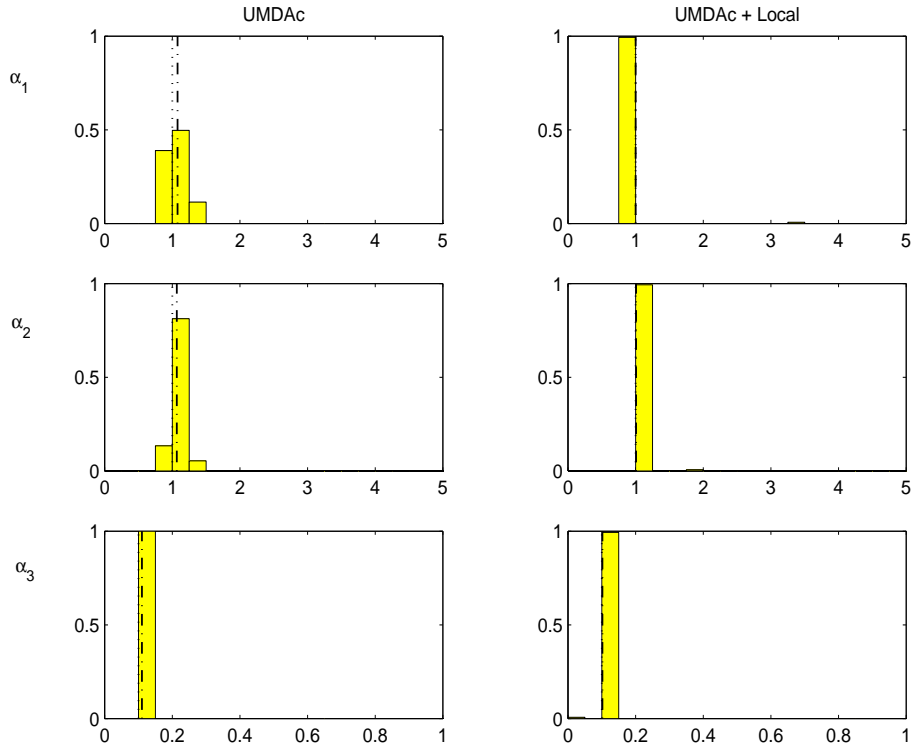}
\caption{\small Histograms of the estimated values for the unknown
parameters $\alpha$ in the model~(\ref{eq:Kon1})-(\ref{eq:Kon2}):
(left) using our proposal for UMDAc, (right) combining the UMDAc and
the local search algorithm.} \label{fig:konno}
\end{figure}

On the contrary, a satisfactory result in the estimation was obtained when
the optimization is carried out via UMDAc with Algorithm~\ref{alg:Global} as it is
shown in Table~\ref{table:UMDAc}. Note that the $100$ estimations of
each parameter are inside of a reduced interval, limited by the lowest and highest estimation values
obtained. The histogram of
the estimator for each parameter is shown at the left column of
Figure~\ref{fig:konno}. In the figure, dot lines correspond to
the real value of the parameters, whereas dash and dot lines
correspond to the mean value of the estimations.

These estimation results can be significantly improved by Algorithm
4 that use each output of Algorithm \ref{alg:Global} as initial value of the
parameters in the local optimization algorithm. Indeed, this is
corroborated in the right column of Figure~\ref{fig:konno} which
shows the histograms of the estimators of $\alpha$ obtained by
Algorithm \ref{alg:GlobalLocal}. It can be observed that the 100 estimations are all grouped very
close to the true value of the parameters.

\subsection[Multiplicative noise example]{Nonlinear model with multiplicative noise}
\label{sec:Cox}

This example is of greater complexity. The state variables of the
diffusion model are strongly dominated by the system noise, the
signal-noise ratio is very low, and the model to be estimated is
over parameterized.

\begin{figure}[htpb]
\centering
\includegraphics[width = 12cm, height = 10cm]{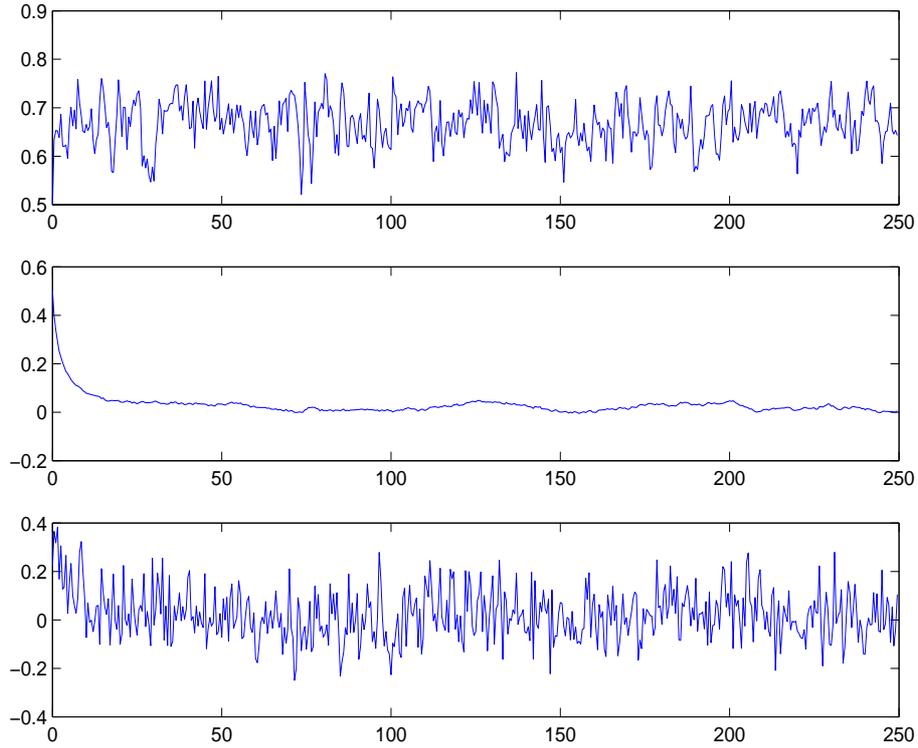}
\caption{ Realization of the solution of the state space model (\ref{eq:cox1})-(\ref{eq:cox2}) Top: $\{x_1(t_j): j=1,\ldots,T\}$,
Center: $\{x_2(t_j): j=1,\ldots,T\}$, Bottom: $\{z_{t_k}:k=0,\ldots ,N\}$}
\label{fig:CoxEq}
\end{figure}

Let us consider the following nonlinear state-space model with multiplicative
noise
\begin{eqnarray}
dx_1 &=&(\alpha_1 +\alpha_2 x_1)dt+ \alpha_3 \sqrt{x_{1}}dw_1 \label{eq:cox1} \\
dx_2 &=&\alpha_4 x_2^2dt+\alpha_5 x_1^2dw_2   \label{eq:cox2}
\end{eqnarray}
\begin{equation}
    z_{t_k}=x_2(t_k)-0.001x_2^3(t_k)+(x_2(t_k)-0.01x_2^2(t_k))\xi_{t_k}+e_{t_k}, \label{eq:CoxObs}
\end{equation}
where $x(t) \in \Re^2$, $z_{t_k}\in \Re$, $\alpha_1 = 1$, $\alpha_2 = -1.5$, $\alpha_3 = 0.1$,
$\alpha_4 = -1$, $\alpha_5 = 0.01$, $x(0) = (x_1(0),x_2(0)) = (0.5,0.5)$, $\{\xi_{t_{k}}:\xi_{t_{k}}\thicksim \mathcal{N}(0,0.01)\}$
and $\{e_{t_{k}}:e_{t_{k}}\thicksim \mathcal{N}(0,0.01)\}$.

In Figure~\ref{fig:CoxEq} a realization $\{x(t_j)=(x_1(t_j),x_2(t_j)): j=1,\ldots,T\}$ of the solution of the
model~(\ref{eq:cox1})-(\ref{eq:cox2}) at instants of time $t_j=jh$, with $h=0.005$ and $T=5\times 10^4$, is shown.
A realization $\{z_{t_k}:k=0,\ldots ,N\}$ of the equation of observations~(\ref{eq:CoxObs}) with
$t_k=k\Delta$, $\Delta =0.5$ and $N=500$ is also shown. It is possible to observe in this figure the strong component of
noise in the signal and in the observations, which gives a high complexity to the estimation problem.

\begin{figure}[htp]
\centering
\includegraphics[width = 15cm, height = 17cm]{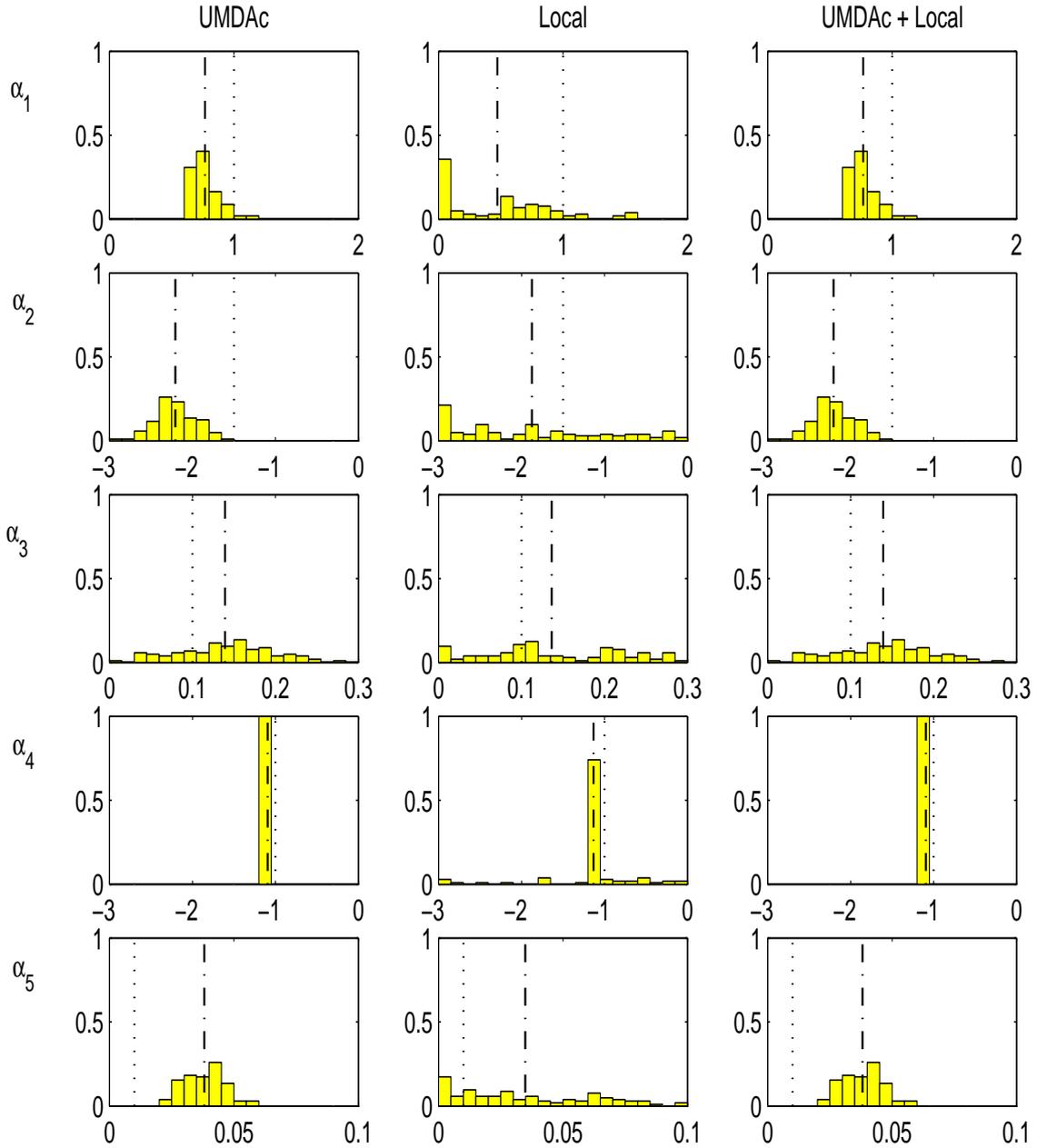}
\caption{ Histograms of the estimated values for the unknown parameters $\alpha$ in the model 
(\ref{eq:cox1})-(\ref{eq:CoxObs}): (left) using our proposal for UMDAc, (center) using the local search algorithm, 
(right) combining the UMDAc and the local search algorithm.}
\label{fig:CoxJK}
\end{figure}

Given the state space model~(\ref{eq:cox1})-(\ref{eq:CoxObs}) and a
single realization of $z_{t_k}$, as is shown in
Figure~\ref{fig:CoxEq}, $100$ estimates of the parameter set $\alpha
= (\alpha_1,\alpha_2,\alpha_3,\alpha_4,\alpha_5)$ were carried out.
Each initial estimate of $\alpha$ was uniformly generated inside of
their definition intervals\[\alpha_1 \in [0,2], \quad \alpha_2 \in
[-3,0], \quad \alpha_3 \in [0,0.3], \quad \alpha_4 \in [-3,0], \quad
\alpha_3 \in [0,0.1].\]

The histograms of the estimated values for $\alpha$ using the local optimization technique of Algorithm \ref{alg:Local} are 
presented in the central column of Figure~\ref{fig:CoxJK}. It is possible to see that the obtained estimation is not useful 
when using a uniformly distributed initial value in the considered intervals for each parameter.

On the other hand, an interesting result in the estimation is obtained with the optimization Algorithms~\ref{alg:Global} 
and \ref{alg:GlobalLocal} via UMDAc, in the sense that all the results are located around certain biased value
of the true parameters. The histograms of the estimated values of $\alpha$ with Algorithm 3 are shown in the column on the left of
Figure~\ref{fig:CoxJK}, while the results obtained with Algorithm \ref{alg:GlobalLocal} are shown in the column on the right. 
In this case, the estimation of Algorithm \ref{alg:GlobalLocal} does not improve the results of Algorithm \ref{alg:Global}, which
confirms that the local optimization algorithm is not useful for this problem.

\begin{figure}[ht]
\centering
\includegraphics[width = 12.5cm, height = 12.5cm]{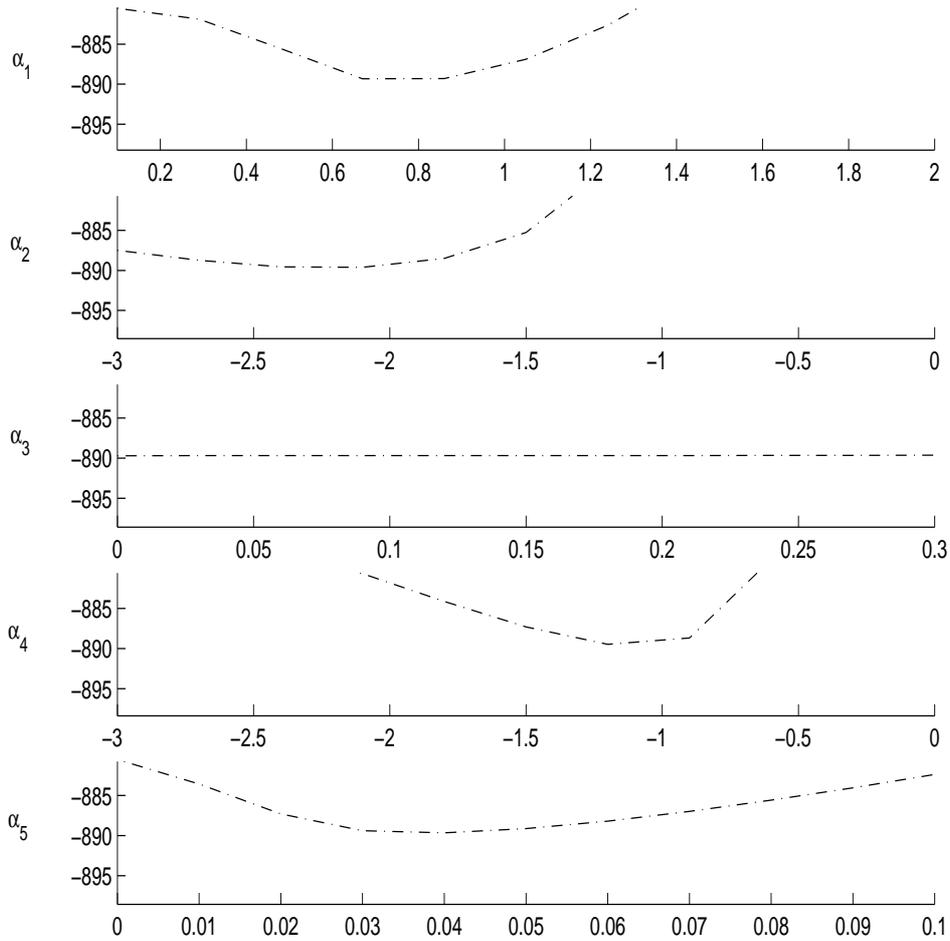}
\caption{\small Evaluation of the fitness function keeping fixed all the coordinates of $\alpha$ except $\alpha_i$
in the model (\ref{eq:cox1})-(\ref{eq:CoxObs}).}
\label{fig:LklhdxParam}
\end{figure}

For a better explanation of the complexity of this optimization problem, the evaluation of the fitness function around the mean
values of the estimated values of the parameters is shown in Figure~\ref{fig:LklhdxParam}, by moving the values of only a
parameter $\alpha_i$. In general, the fitness function is almost flat, which is even more evident for the third coordinate. The
existence of a local minimum for the third parameter could be observed making a bigger rescaling in the considered interval.
However, note that in this extreme situation the mean value of the estimated parameters $\alpha_1$, $\alpha_2$, $\alpha_4$, and
$\alpha_5$ are quite close to the minimum value of the fitness function for each parameter, which reveals the good performance of
the Algorithms \ref{alg:Global} and \ref{alg:GlobalLocal}. This demonstrates the robustness of the UMDAc
in the case of having small deviations from the typical assumptions required for its use.

\section[Conclusions]{Conclusions} \label{sec:CONCLU}

In this paper, we have considered two optimization methods based on
the Estimation of Distribution Algorithms for computing the
Innovation Estimators of the unknown parameters of diffusion
processes given a set of discrete and noisy observations. The first
method is exclusively based on a variant of the known Univariate
Marginal Distribution Algorithm in continuous domain, whereas the
second method includes a refinement for the outputs of the first one
via a local optimization algorithm. The performance of these two
optimization methods were evaluated in the parameter estimation of
two types of diffusion models with complex nonlinear and stochastic
dynamics. The numerical simulations demonstrate the feasibility of
the considered method for the parameter estimation in situations
where local optimization algorithms fail. This is particularly
relevant in practice when adequate initial values for the parameters
to be estimated are not available or when the diffusion model has
highly nonlinear parameters to be estimated from highly noisy
observations.

While our results show that UMDAc is able to deal with optimization
scenarios where the commonly applied local optimization methods
fail, there are still room for improvement. In particular, other
EDAs that explicitly model and exploit multivariate interactions
between the parameters of the fitness function are worth to be
evaluated in the computation of Innovation Estimators of
dif\mbox{}fusion processes. We leave this question as a line of
future research.

\begin{acknowledgements}
Z.G.A. acknowledges partial financial support from the \emph{Funda\c c\~ao Carlos Chagas Filho de Amparo \`a
Pesquisa do Estado do Rio de Janeiro} (FAPERJ).
\end{acknowledgements}


\end{document}